\documentclass[reqno,11pt]{amsart}
\usepackage{amsmath, latexsym, amsfonts, amssymb, amsthm, amscd}
\usepackage{mathrsfs,enumerate}
\usepackage{amsmath,amssymb}

\numberwithin{equation}{section}

\newtheorem{theorem}{Theorem}[section]
\newtheorem{lemma}[theorem]{Lemma}
\newtheorem{prop}[theorem]{Proposition}

\newcommand{\Probabilities}[1]{\mathscr P(#1)}
\newcommand{\ProbabilitiesTwo}[1]{\mathscr P_2(#1)}
\newcommand{\RelativeEntropy}[2]{\mathcal H(#1|#2)}

\newcommand{\gga}{\gamma}            

\newcommand{\gep}{\varepsilon}       

\newcommand{\cE}{{\ensuremath{\mathcal E}} }
\newcommand{\cH}{{\ensuremath{\mathcal H}} }

\newcommand{\cN}{{\ensuremath{\mathcal N}} }

\newcommand{\cL}{{\ensuremath{\mathcal L}} }

\newcommand{\E}{{\ensuremath{\mathbb E}} }

\newcommand{\N}{{\ensuremath{\mathbb N}} }

\newcommand{\bbP}{{\ensuremath{\mathbb P}} }
\newcommand{\Q}{{\ensuremath{\mathbb Q}} }
\newcommand{\R}{{\ensuremath{\mathbb R}} }

\newcommand{\V}{{\ensuremath{\mathbb V}} }

\author{Lorenzo Zambotti}
\address{Laboratoire de Probabilit{\'e}s et Mod\`eles Al\'eatoires (CNRS U.M.R. 7599) \\ Universit{\'e} Paris 6
-- Pierre et Marie Curie, U.F.R. Mathematiques, Case 188, 4 place
Jussieu, 75252 Paris cedex 05, France }
\email{zambotti\@@ccr.jussieu.fr}

\date{}

\title{Fluctuations for a conservative interface model on a wall}

\begin{document}

\begin{abstract}
We consider an effective interface model on a hard wall in (1+1)
dimensions, with conservation of the area between the interface
and the wall. We prove that the equilibrium fluctuations of the height variable
converge in law to the solution of a SPDE with reflection and
conservation of the space average. The proof is based on recent
results obtained with L. Ambrosio and G. Savar\'e on stability
properties of Markov processes with log-concave invariant measures.
\\
\\
2000 \textit{Mathematics Subject Classification: 60K35; 60H15;
82B05}
\end{abstract}

\keywords{Equilibrium fluctuations; Interface model; Stochastic
partial differential equations; hard wall}

\maketitle

\section{Introduction}

This paper concerns fluctuations of a $\nabla\phi$ interface model
on a hard wall with conservation of the area between the interface
and the wall. The system is defined on the one-dimensional lattice
$\Gamma_N:=\{1,2,\ldots,N\}$ and the location of the interface at
time $t$ is represented by the height variables
$\phi_t=\{\phi_t(x),x\in\Gamma_N\}\in
\Omega_N^+:=[0,\infty)^{\Gamma_N}$ measured from the wall
$\Gamma_N$.

In order to describe the dynamics of $\phi_t$ we need some notation.
Let $\{(w_t(x))_{t\geq 0}:x=1,\ldots,N\}$ be independent standard
Brownian motions and define the $N\times N$ matrices
\[
\sigma \, := \, \begin{pmatrix}
-1 & \  &   &    &\\
\ \  1  & \ \cdot &  &    &\\
   &  \ \cdot  & \ \cdot &  &\\
   &    & 1  &  -1 &  \\
   &   &   &   \ \ 1  & \ 0
\end{pmatrix}, \qquad
\sigma^T \, := \, \begin{pmatrix}
-1 & \ 1 &   &    &\\
   & \ \cdot & \ \cdot &    &\\
   &    & \ \cdot & \ \cdot &\\
   &    &   &  -1 & \ 1 \\
   &   &   &   \ 0  & \ 0
\end{pmatrix}
\]
Then the dynamics of $(\phi_t(x):x\in\Gamma_N)_{t\geq 0}$, height
from the wall of the reflected interface, is governed by the
stochastic differential equation of the Skorohod type
\begin{equation}\label{sde}
d\phi_t \, = \, - \, \sigma\sigma^T \{\sigma V'(\sigma^T\phi_t)\, dt
+ \, dl_t \} + \, \sqrt 2 \,\sigma \, dw_t
\end{equation}
for all $x\in \Gamma_N$, subject to the conditions
\begin{equation}\label{constraint}
\begin{split}
 \phi_t(x) \, \geq \, 0, \qquad & t\mapsto l_t(x) \ {\rm continuous \
and \ non-decreasing}, \qquad l_0(x)=0,
\\ &
\int_0^\infty \phi_t(x) \, dl_t(x) \, = \, 0, \qquad x\in
\Gamma_N.
\end{split}
\end{equation}
We refer to \cite{fu} for an introduction to interface models.

\smallskip\noindent
Throughout the paper the potential $V$ satisfies the following conditions
\begin{itemize}
\item[(V1)] (convexity) $V\in C^2(\R)$ is convex and
$\lim_{|r|\to\infty} V(r)=+\infty$.
\end{itemize}
Notice that for a convex $V$
\[
\lim_{|r|\to\infty} V(r)=+\infty \, \Longleftrightarrow \, \int
\exp(-V)\, dr<\infty \, \Longleftrightarrow \, V(r)\geq a+b|r| \quad
\forall \ r\in\R,
\]
for some $a\in\R$ and $b>0$. In particular we have
\begin{equation}\label{qdet}
q \, := \, \int_{\mathbb R} r^2 \, \exp(-V(r)) \, dr \, < \, \infty
\end{equation}
\begin{itemize}
\item[(V2)] (normalization), \quad
${\displaystyle \int_{\mathbb R} \exp(-V(r)) \, dr \, = \, 1}$.
\item[(V3)] ($0$ mean), \quad
${\displaystyle \int_{\mathbb R} r \, \exp(-V(r)) \, dr \, = \, 0}$.
\end{itemize}
The normalization (V2) does not affect equation
(\ref{sde}), where only $V'$ appears.

\medskip\noindent We shall prove in the following sections existence
and uniqueness of solutions of \eqref{sde} and other properties.

\subsection{The main result}
For any $N\in{\mathbb N}$ we set $\Lambda_N:{\mathbb R}^N\mapsto
L^2(0,1)$,
\begin{equation}\label{rescaling}
\Lambda_N(\phi)(\theta) \, := \, \frac 1{\sqrt N} \, \phi(\lfloor
N\theta\rfloor +1), \quad \theta\in[0,1),
\end{equation}
where $\lfloor\cdot\rfloor$ denotes the integer part, and we define
the spaces
\[
H_N=\Lambda_N({\mathbb R}^N)\subset L^2(0,1), \qquad \Omega^+_N :=
({\mathbb R}_+)^N, \qquad  K_N:=\Lambda_N(\Omega^+_N).
\]
Notice that $K_N$ can be identified with the space of non-negative
functions on $[0,1)$ being constant on $I(x)=[(x-1)/N,x/N)$ for all
$x\in\Gamma_N$.

For all $k\in K_N$ and $t\geq 0$ we define now the rescaled
interface $\Phi^N$
\[
\Phi^N_t \, := \, \Lambda_N\left( \phi_{N^4t} \right),
\qquad \Phi^N_0 := \Lambda_N\left( \phi_0 \right).
\]
In other words
\[
\Phi^N_t(\theta) \, = \, \frac 1{\sqrt N} \, \phi_{N^4t}(\lfloor
N\theta\rfloor +1), \quad \theta\in[0,1).
\]
In the main result of this paper, i.e. Theorem \ref{main} below, we
state the weak convergence of $\Phi^N$ to the unique solution $u$ of
the following stochastic Cahn-Hilliard equation on $[0,1]$ with
homogeneous Neumann boundary condition and reflection at $u=0$
\begin{equation}\label{0}
\left\{ \begin{array}{ll} {\displaystyle \frac{\partial u}{\partial
t}=-\frac{\partial^2 }{\partial \theta^2} \left(\frac 1q\,
\frac{\partial^2 u}{\partial \theta^2} +\eta \right)  + \sqrt 2 \,
\frac{\partial}{\partial\theta} \dot{W}, }
\\ \\
{\displaystyle \frac{\partial u}{\partial
\theta}(t,0)=\frac{\partial u}{\partial \theta}(t,1)=
\frac{\partial^3 u}{\partial \theta^3}(t,0)=\frac{\partial^3
u}{\partial \theta^3}(t,1)=0 ,}
\\ \\
u(0,\theta)=u_0(\theta), \quad \theta\in[0,1]
\end{array} \right.
\end{equation}
where $\dot W$ is a space-time white noise on
$[0,+\infty)\times[0,1]$, $u$ is a continuous function of
$(t,\theta)\in [0,+\infty)\times[0,1]$, $\eta$ is a locally finite
positive measure on $(0,+\infty)\times[0,1]$, subject to the
constraint
\begin{equation}\label{contact}
u\geq 0, \qquad \int_{(0,+\infty)\times[0,1]} u \, d\eta \, = \, 0.
\end{equation}
Such equation has been studied in \cite{deza}, see Proposition \ref{deza}
below.

With an abuse of notation, we say that a sequence of measures $({\bf
P}_n)$ on $C([a,b];L^2(0,1))$ converges weakly in
$C([a,b];L^2_w(0,1))$ if, for all $m\in\N$ and $h_1,\ldots,h_m\in
C^1([0,1])$, the process $(\langle X_\cdot,h_i\rangle_{L^2(0,1)}, \,
i=1,\ldots,m)$ under $({\bf P}_n)$ converges weakly in
$C([a,b];\R^m)$ as $n\to\infty$.

Then we can state the main result of this paper.
\begin{theorem}\label{main}
If $\Phi^N_0\to u_0$ in $L^2(0,1)$ as $N\to\infty$ with
\[
\Phi^N_0 \geq 0,  \qquad
\int_0^1\Phi^N_0(\theta) \, d\theta=c>0 \qquad \forall\ N\in\N,
\]
then, for any $0<\varepsilon\leq T<\infty$, the law of  $(\Phi^N_t,
t\in[\varepsilon,T])$ converges to the
law of the unique solution $u$ of \eqref{0}, weakly in $C([\varepsilon,T];L^2_w(0,1))$.
\end{theorem}

\subsection{A conservative dynamics}

The starting point of this work is the paper by Funaki and Olla
\cite{fuol}. In that paper, the following $\nabla\phi$ interface
model on a hard wall is considered
\begin{equation}\label{sdesko}
d\phi_t(x)  =  - \, \sigma V'(\sigma^T\phi_t) \, dt + \, dl_t(x) \,
+ \, \sqrt 2 \, dw_t(x), \qquad x\in \Gamma_N,
\end{equation}
with constraints analogous to \eqref{constraint} and Dirichlet
boundary condition $\phi_t(0)=\phi_t(N+1)=0$. Using the definition
\eqref{rescaling}, it is then proven that in the stationary case,
the process $(\Lambda_N(\phi_{N^2t}), t\geq 0)$ converges in law as
$N\to\infty$ to the law of the unique stationary solution of the
second order equation
\begin{equation}\label{alpha0}
\left\{ \begin{array}{ll} {\displaystyle \frac{\partial u}{\partial
t}=\frac 1q \frac{\partial^2 u}{\partial \theta^2} + \eta
 + \sqrt 2\,  \frac{\partial^2 W}{\partial t\partial\theta} }
\\ \\
u(t,0)=u(t,1)=0, \quad t\geq 0
\\ \\
u\geq 0, \ d\eta\geq 0, \ \int u\, d\eta=0
\end{array} \right.
\end{equation}

At the end of the introduction of \cite{fuol}, it is remarked that
it would be more natural to consider a stochastic dynamics
conserving the area between the interface and the wall, namely
$\sum_x \phi(x)$. Such conservative dynamics, but without the hard
wall constraint, has indeed been studied in \cite{nishi} and
\cite{nishi2}, where respectively hydrodynamic limit and large
deviations are considered; the hydrodynamic scaling limit of the
interface is the solution of a fourth-order equation, as predicted
in \cite{spohn}.

The SDE \eqref{sde} combines the hard wall and the conservation of
volume  constraints; indeed, $\sigma^T{\bf 1}=0$, where ${\bf
1}=(1,\ldots,1)\in\R^N$, and it is easy to see that
\[
d\left[\sum_{x=1}^N \phi_t(x) \right] =  \sum_{x=1}^N
\left[\sigma^T{\bf 1}\right](x) \left\{-\left[\sigma^T \{\sigma
V'(\sigma^T\phi_t)\, dt + \, dl_t \}\right](x) + \, \sqrt 2 \,
dw_t(x) \right\} = 0.
\]

The main novelty of this paper is the use of a technique recently
developed in \cite{asz} for the convergence in law of stochastic
processes associated with symmetric Dirichlet forms of gradient type
and with log-concave invariant measures; see section \ref{asz}
below. The general principle is in fact very simple: this class
of reversible dynamics is parametrized by two objects, the invariant
measure and the scalar product of the Hilbert space which defines
the gradient. If such objects converge (in a sense te be made precise),
it is natural to conjecture
that the associated processes converge; the results of \cite{asz}
confirm this conjecture in the case of log-concave reference measures:
see section \ref{asz} below.

The solutions of equations \eqref{sde}, \eqref{0},
\eqref{sdesko} and \eqref{alpha0} are all in this class and the
techniques of \cite{asz} give a general framework to prove results
like Theorem \ref{main} or the convergence result of \cite{fuol}. We
recall that \cite{fuol} is based on monotonicity properties, which
are rather special properties of \eqref{sdesko}-\eqref{alpha0}, not
shared by \eqref{sde}-\eqref{0}. One can notice that, given the
general results of \cite{asz}, the proof of convergence of equilibrium
fluctuations as in \cite{fuol} and in this paper becomes much
easier.

We also notice that Theorem \ref{main} is comparatively stronger
than the analogous statement in \cite{fuol}. Indeed, we consider a
convex microscopic interaction potential $V$, instead of a strictly
convex and symmetric one. Moreover the convergence is proven not only in the
stationary case, but for any sequence of initial conditions which
converge under the rescaling \eqref{rescaling}. Using the techniques
of this paper, one could improve correspondingly the results of
\cite{fuol}.

Finally, we notice that the boundary conditions we consider are of
Neumann type, like in \cite{cdg}, while many other papers
consider the Dirichlet (see e.g. \cite{fuol}) or the periodic
(see e.g. \cite{nishi}) case. The case of periodic boundary condition
could be proven with no additional difficulty with the techniques of this paper.
Indeed, like in the Neumann case, the invariant measure of the limit
SPDE is absolutely continuous w.r.t. the Gaussian invariant probability
measure of the linear SPDE (i.e. without reflection). The weak convergence
of the rescaled stationary measures is then a simple
consequence of a standard invariance principle: see the proof of Proposition
\ref{asz2}.

For the case of homogeneous Dirichlet boundary conditions, on the contrary, the
invariant measure of the limit SPDE is singular w.r.t. the Gaussian invariant probability
measure of the linear SPDE, due to the interplay of the homogeneous boundary conditions
and the non-negativity constraint. This makes the convergence of the
rescaled invariant measures more delicate. In fact,
we could prove the results of this paper for
Dirichlet boundary condition, if we could prove the following
invariance principle: we consider
a random walk $S_n=X_1+\cdots+X_n$, $n=1,\ldots,N$, with step distribution $X_i\sim
e^{-V}dx$, conditioned to be non-negative (i.e. $S_1,\ldots,S_N\geq 0$),
to be $0$ at time $N$ (i.e. $S_N=0$) and to have a fixed sum
(i.e. $\sum_{n=1}^N S_n=c N^{3/2}$, $c>0$); then we would like to prove that
such processes converge under Brownian rescaling as $N\to\infty$ to
a Brownian excursion $e$ conditioned to have integral $c$
(i.e. $\int_0^1 e_x \, dx =c$). Since we have not found a proof for this
invariance principle, we restrict to the Neumann case, for which we
can prove convergence of the stationary measures. In the Dirichlet boundary condition case the
limit SPDE would be an analog of \eqref{0}, with boundary
conditions
\[
\displaystyle u(t,0)=u(t,1)= \frac{\partial^3 u}{\partial
\theta^3}(t,0)=\frac{\partial^3 u}{\partial \theta^3}(t,1)=0 ,
\]
i.e. Dirichlet for $u$ and Neumann for
$\frac{\partial^2u}{\partial\theta^2}$. Such equation is studied in
\cite{z}.

\section{A general convergence result}\label{asz}

In this section we recall the results of \cite{asz}, already
mentioned in the introduction. It turns out that the processes
$(\phi_t)$ and $(u(t,\cdot)$, solutions of \eqref{sde} and \eqref{0}
respectively, are both {\it monotone gradient systems}, i.e. the
equation they satisfy can be interpreted as follows
\[
dX = -\nabla U(X) \, dt + \sqrt 2 \, dW
\]
where $W$ is a Wiener process in a Hilbert space $H$ and $U:H\mapsto
\R\cup\{+\infty\}$ is a convex potential. These processes are
reversible and associated with a gradient-type Dirichlet form. The
general results of existence and convergence of such processes given
in \cite{asz}, have a nice application in the present setting. Hence
we devote this section to recall them.

\medskip Let $H$ be a separable Hilbert space with scalar product
$\langle\cdot,\cdot\rangle_H$ and let $\gamma$ be a probability measure on
$H$. We suppose that $\gamma$ is {\it log-concave}, i.e. for all
pairs of open sets $B,\, C\subset H$
\begin{equation}\label{deflogconc}
\log\gamma\left((1-t)B+tC\right)\geq
(1-t)\log\gamma(B)+t\log\gamma(C) \qquad\forall t\in (0,1).
\end{equation}
If $H=\R^k$, then the class of log-concave probability measures
contains all measures of the form (here $\cL_k$ stands for Lebesgue
measure)
\begin{equation}
  \label{eq:basic_example}
  \gamma :=  \frac 1Z \, e^{-U} \cL_k,
\end{equation}
where $U:H=\R^k\to\R\cup\{+\infty\}$ is convex and $Z:= \int_{\R^k}
e^{-U} \, dx<+\infty$, see Theorem~9.4.11 in \cite{ags}, in
particular all Gaussian measures. Notice that the class of
log-concave measures is closed under weak convergence. Moreover, if
$\gamma$ is log-concave and $K$ is a convex set with $\gamma(K)>0$,
then the conditional measure $\gamma(\cdot | K):=\gamma(\cdot\cap
K)/\gamma(K)$ is also log-concave.

We denote the support of $\gamma$ by $K=K(\gamma)$ and the smallest
closed affine subspace of $H$ containing $K$ by $A=A(\gamma)$. We
write canonically
\begin{equation}\label{defah}
A \, = \, H^0 \, + \, h^0,
\end{equation}
where $H^0=H^0(\gamma)$ is a closed linear subspace of $H$ and
$h^0=h^0(\gamma)$ is the element of minimal norm in $A$. We endow
$H^0$ with the scalar product $\langle \cdot,\cdot\rangle_{H^0}$
induced by $H$.

We want to consider a stochastic processes with values in $A$ and
reversible with respect to $\gamma$. We denote by $C_b(H)$ the space
of bounded continuous functions in $H$ and by $C_b^1(A)$ the space
of all $\Phi:A\mapsto\R$ which are bounded, continuous and Fr\'echet
differentiable. To $\varphi\in C_b^1(A)$ we associate a gradient
$\nabla_{H^0}\varphi:A\mapsto H^0$, defined by
\begin{equation}\label{gradiente}
\left.\frac d{d\varepsilon} \, \varphi(k+\gep \, h) \right|_{\gep=0}
\, = \, \langle\nabla_{H^0}\varphi(k),h\rangle_{H^0}, \qquad \forall
\ k\in A, \ h\in H^0.
\end{equation}

We denote by $X_t:K^{[0,+\infty[}\to K$ the coordinate process
$X_t(\omega):=\omega_t$, $t\geq 0$. Finally, we denote the set of
probability measures on $H$ by $\Probabilities{H}$ and we set
$$
\ProbabilitiesTwo{H}:=\left\{\mu\in\Probabilities{H}:\
\int_H\Vert x\Vert^2_H\,d\mu(x)<\infty\right\},
$$Then we recall one of the main
results of \cite{asz}.
\begin{theorem}[Markov process and Dirichlet form associated with
$\gamma$ and $\|\cdot\|_{H^0}$]\label{main1}
$ $
\begin{itemize}
\item[(a)]
The bilinear form ${\cE}={\cE}_{\gamma,\|\cdot\|_{H^0}}$ given by
\begin{equation}\label{diri}
{\cE}(u,v) \, := \, \int_{K} \langle\nabla_{H^0} u,\nabla_{H^0}
v\rangle_{H^0} \, d\gamma, \qquad u,\,v\in C^1_b(A),
\end{equation}
is closable in $L^2(\gamma)$ and its closure $(\cE,D(\cE))$ is a
symmetric Dirichlet Form. Furthermore, the associated semigroup
$(P_t)_{t\geq 0}$ in $L^2(\gamma)$ maps $L^\infty(\gamma)$ in
$C_b(K)$. \smallskip
\item[(b)] There exists a unique Markov family $(\bbP_x:x\in K)$ of probability
measures on $K^{[0,+\infty[}$ associated with $\cE$. More precisely,
$\E_x[f(X_t)]=P_tf(x)$ for all bounded Borel functions and all $x\in
K$. \smallskip
\item[(c)] For all $x\in K$, $\bbP_x^*\left(C(]0,+\infty[;H)\right)=1$ and
$\E_x[\|X_t-x\|^2]\to 0$ as $t\downarrow 0$. Moreover,
$\bbP_x^*\left(C([0,+\infty[;H)\right)=1$ for $\gamma$-a.e. $x\in
K$. \smallskip
\item[(d)] $(\bbP_x:x\in K)$ is reversible with respect to $\gamma$,
i.e. the transition semigroup $(P_t)_{t\geq 0}$ is symmetric in
$L^2(\gamma)$; moreover $\gamma$ is invariant for $(P_t)$, i.e.
$\gamma(P_tf)=\gamma(f)$ for all $f\in C_b(K)$ and $t\geq 0$.
\item[(e)] If $\gamma\in\ProbabilitiesTwo{H}$, then $\gamma$
is the only invariant probability measure for $(P_t)$ in
$\ProbabilitiesTwo{H}$.
\end{itemize}
\end{theorem}
We shall see below that the solutions of \eqref{sde}, \eqref{0},
\eqref{sdesko} and \eqref{alpha0} are all particular cases of the
class of Markov processes described in Theorem \ref{main1}. This
fact will be crucial in the proof of Theorem \ref{main}.

\medskip\noindent
We consider now a sequence $(\gamma_N)$ of log-concave probability
measures on $H$ such that $\gamma_N$ converge weakly in $H$ to
$\gamma$. We denote by $K_N$ the support of $\gamma_N$, and by $A_N$
the smallest closed affine subspace of $H$ containing $K_N$. We
suppose that $A_N\subseteq A$ for all $N$.

We write $A_N=h^0_N+H_N^0$, where $h^0_N\in A_N$ and $H_N^0\subseteq
H^0$ is a closed linear subspace of $H$. We want to consider
situations where each $H_N^0$ is a Hilbert space endowed with a
scalar product $\langle\cdot,\cdot\rangle_{H_N^0}$, possibly
different from the scalar product induced by $H^0$. In order to
ensure that this family of scalar products converges (in a suitable
sense) to the scalar product of $H^0$ as $N\to\infty$, we will make
the following assumptions.
\begin{enumerate}
\item There exists a finite constant $\kappa\geq 1$ such that
\begin{equation}\label{bka0}
\frac{1}{\kappa}\Vert h\Vert_{H^0}\leq \Vert
h\Vert_{H_N^0}\leq\kappa\Vert h\Vert_{H^0} \qquad\forall \ h\in
H_N^0, \ N \in\N.
\end{equation}
\item Denoting by $\Pi_N:H^0\to H_N^0$ the orthogonal
projections induced by the scalar product of $H^0$, we have
\begin{equation}\label{sig0}
\lim_{N\to\infty}\Vert\Pi_N h\Vert_{H_N^0}=\Vert h\Vert_{H^0}
\qquad\forall \ h\in H^0.
\end{equation}
\end{enumerate}
These assumptions guarantee in some weak sense that the geometry of
$H_N^0$ converges to the geometry of $H^0$; the case when all scalar
products coincide with $\langle\cdot,\cdot\rangle_H$, $H_N^0\subset
H_{N+1}^0$ and $\cup_N H_N^0$ is dense in $H^0$ is obviously
included.

Let $(\bbP_x^N:x\in K_N)$ (respectively $(\bbP_x:x\in K)$) be the
Markov process in $[0,+\infty[^{K_N}$ associated to $\gamma_N$
(resp. in $[0,+\infty[^K$ associated to $\gamma$) given by
Theorem~\ref{main1}. We denote by $\bbP_{\gamma_N}^N:=\int
\bbP_{x}^N \, d\gamma_N(x)$ (resp. $\bbP_\gamma:= \int \bbP_{x} \,
d\gamma(x)$) the associated stationary measures.

With an abuse of notation, we say that a sequence of measures $({\bf
P}_n)$ on $C([a,b];H)$ converges weakly in $C([a,b];H_w)$ if, for
all $m\in\N$ and $h_1,\ldots,h_m\in H$, the process $(\langle
X_\cdot,h_i\rangle_H, \, i=1,\ldots,m)$ under $({\bf P}_n)$
converges weakly in $C([a,b];\R^m)$ as $n\to\infty$.

In this setting we have the following stability and tightness
result, also proven in \cite{asz}.
\begin{theorem}[Stability and tightness]\label{main3}
Suppose that $\gamma_N\to\gamma$ weakly in $H$ and that the norms of
$H_N^0$ satisfy \eqref{bka0} and \eqref{sig0}. Then, for any $x_N\in
K_N$ such that $x_N\to x\in K$ in $H$, for any $0<\varepsilon\leq
T<+\infty$, $\bbP_{x_N}^N\to\bbP_x$ weakly in
$C([\varepsilon,T];H_w)$;
\end{theorem}
This stability property means that the weak convergence of the
invariant measures $\gamma_N$ and a suitable convergence of the
norms $\|\cdot\|_{H_N^0}$ to $\|\cdot\|_{H^0}$ imply the convergence
in law of the associated processes, starting from any initial
condition.

\medskip\noindent
We recall that the above results, proven in \cite{asz}, are based on
the interpretation of the Markov semigroup $(P_t)$ as the solution
of a gradient flow in $\ProbabilitiesTwo{H}$ with respect to the
relative entropy functional $\RelativeEntropy{\cdot}{\gamma}$ in the
Wasserstein metric: see  \cite{asz} for details.

In the rest of the paper we show how the results of this
section apply to Theorem \ref{main}.

\section{The microscopic dynamics}

On $\R^N$ we consider the canonical scalar product and we denote it
by $\langle\cdot,\cdot\rangle_{\R^N}$, with associated norm
$\|\cdot\|_{\R^N}$.

We define ${\bf 1}:=(1,\ldots,1)\in\R^N$ and the vector space
$\V_N:=\{v\in\R^N: v_1+\cdots+v_N=0\} = {\bf 1}^\perp$. It is easy
to see that the kernels of $\sigma$ and $\sigma^T$ are respectively
$\text{Ker}(\sigma)=\{(0,\ldots,0,t): t\in\R\}$ and
$\text{Ker}(\sigma^T)= \{t\cdot {\bf 1}\in\R^N: t\in\R\}$; it
follows that the image of $\sigma$ is
$\text{Im}(\sigma)=(\text{Ker}(\sigma^T))^\perp= \V_N$ and that
$\text{Ker}(\sigma)\cap\V_N=\{0\}$; therefore
$\sigma:\V_N\mapsto\V_N$ is bijective, $\sigma^{-1}:\V_N\mapsto\V_N$
is well defined and we can define the scalar product in $\V_N$
\[
\langle v_1,v_2 \rangle_{\V_N} \, := \, \langle\sigma^{-1}
v_1,\sigma^{-1} v_2\rangle_{\R^N}, \qquad \forall \ v_1,v_2\in \V_N.
\]
We want now to give a useful representation of $\langle \cdot,\cdot
\rangle_{\V_N}$. Let $(B_t, t\geq 0)$ be a standard Brownian motion
and set
\begin{equation}\label{D}
D_i := B_i -\frac{B_1+B_2+\cdots+B_N}N, \quad i=1,\ldots,N, \qquad
D:=(D_1,\ldots,D_N)\in\V_N.
\end{equation}
\begin{lemma}\label{repre}
For all $v\in\V_N$
\[
\|v\|^2_{\V_N} \, = \, \E\left[ \langle v,D\rangle_{\R^N}^2 \right]
\, = \, \sum_{i=1}^{N-1} \left( \sum_{j=1}^i v_j \right)^2.
\]
\end{lemma}
\noindent{\it Proof}. Let $V\in\V_N$ such that $\sigma V=v$. Then
$\|v\|^2_{\V_N} \, = \, \|V\|^2_{\R^N}$. Moreover $V_i=\sum_{j=1}^i
v_j$, $i=1,\ldots,N$, and in particular $V_N=0$ since $v\in\V_N$.
Since $\sigma^T D=(B_2-B_1,\ldots,B_N-B_{N-1},0)$ and $V_N=0$
\[
\E\left[ \langle v,D\rangle_{\R^N}^2 \right] \, = \, \E\left[
\langle V,\sigma^TD\rangle_{\R^N}^2 \right] = \|V\|^2_{\R^N} =
\|v\|^2_{\V_N}. \qed
\]

\medskip\noindent Recall that $\{(w_t(x))_{t\geq 0}: x=1,\ldots,N\}$ is an independent
family of standard Brownian motions; then $w=(w(1),\ldots,w(N))$ is
a Wiener process in $\R^N$ and $\sigma w$ is a Wiener process in
$\V_N$, i.e. for all $t\geq 0$
\[
\E\left[ \langle h,w_t\rangle_{\R^N}^2 \right] = t\|h\|^2_{\R^N},
\quad \forall \ h\in\R^N, \qquad \E\left[ \langle v,\sigma
w_t\rangle_{\V_N}^2 \right] = t\|v\|^2_{\V_N}, \quad \forall \
v\in\V_N. 
\]

\begin{lemma}\label{exun}
For all $\phi_0\in K_N$ there exists a unique pair
$(\phi_t,l_t)_{t\geq 0}$, solution of {\rm (\ref{sde})}. We use the
notation $\phi(t,\phi_0)=\phi_t$, $t\geq 0$.
\end{lemma}
\noindent{\it Proof}. We start by (pathwise) uniqueness. Let $(\phi,
l)$ and $(\overline{\phi},\overline{l})$ be solutions of (\ref{sde})
with initial condition $\phi_0$, resp. $\overline{\phi}_0$. Setting
$\psi_t:=\phi_t-\overline{\phi}_t$, by It\^o's formula we obtain
\[
d\langle \psi_t,{\bf 1} \rangle_{\R^N} = \langle \sigma^T{\bf 1},
-\, \sigma^T \{\sigma
(V'(\sigma^T\phi_t)-V'(\sigma^T\overline\phi_t))\, dt + \,
dl_t-d\overline{l}_t \}\rangle_{\R^N} = 0
\]
so that $\langle \psi_t,{\bf 1} \rangle=0$ for all $t\geq 0$ and
therefore $\psi_t\in\V_N$. Then, again by It\^o's formula
\begin{eqnarray*} & &
d\langle\psi_t,\psi_t\rangle_{\V_N}  \, = \, - \langle
\sigma^T\psi_t, V'(\sigma^T\phi_t)-V'(\sigma^T\overline{\phi}_t)
\rangle \, dt + \langle \psi, dl_t-d\overline{l}_t \rangle_{\R^N} \,
\leq \, 0
\end{eqnarray*}
since $V'$ is monotone non-decreasing and by \eqref{constraint}.

For or existence of (strong) solutions, we can refer to \cite{cepa}.
Indeed, setting ${\bf 1}_\phi:=\langle\phi_0,{\bf 1}\rangle_{\R^N}
{\bf 1}$ and
$\zeta_t:=\phi_t-{\bf 1}_\phi$, \eqref{sde} is equivalent to
\begin{equation}\label{sdeze}
d\zeta_t \, = \, - \, \sigma\sigma^T \{\sigma
V'(\sigma^T(\zeta_t+{\bf 1}_\phi))\,
dt + \, dl_t \} + \, \sqrt 2 \,\sigma \, dw_t
\end{equation}
for all $x\in \Gamma_N$, subject to the conditions
\[
\begin{split}
 \zeta_t(x)+\langle\phi_0,{\bf 1}\rangle_{\R^N} \, \geq \, 0, \qquad & t\mapsto l_t(x) \ {\rm continuous \
and \ non-decreasing}, \qquad l_0(x)=0,
\\ &
\int_0^\infty \left(\zeta_t(x)+\langle\phi_0,{\bf
1}\rangle_{\R^N}\right) \, dl_t(x) \, = \, 0, \qquad x\in \Gamma_N.
\end{split}
\]
Equation \eqref{sdeze} is a Skorohod problem in the convex set
$[0,\infty[^{\Gamma_N}\cap \V_N$; in other words, $\zeta$ solves the
stochastic differential {\it inclusion}
\[
d\zeta \in -\partial U(\zeta_t)\, dt + \sqrt 2 \,\sigma \, dw_t
\]
where $U:\V_N\mapsto\R$ is the convex potential
\[
U(\zeta) := \left\{ \begin{array}{ll} {\displaystyle \sum_{x=2}^N
V(\zeta(x)-\zeta(x-1)), \qquad \text{if } \ \zeta+{\bf 1}_\phi \in [0,\infty[^{\Gamma_N}\cap \V_N }
\\ \\
\qquad \qquad    +\infty, \qquad \qquad \qquad \qquad \qquad
\text{otherwise,}
\end{array} \right.
\]
see in particular Proposition 3.1 in \cite{cepa}. Therefore
existence of a strong solution of \ref{sdeze} follows from Theorem
5.1 of \cite{cepa}. \qed

\section{The microscopic invariant measure}

In this section we study invariant measures of \eqref{sde} and
the associated Dirichlet forms. Since \eqref{sde} conserves the sum
$\sum_{x=1}^N \phi_t(x)=\sum_{x=1}^N \phi_0(x)$ for all $t\geq 0$,
each subspace $\V_N^c=\V_N+c{\bf 1}$, with $c>0$,
supports an invariant measure. Therefore it is natural to fix $c>0$
and consider only initial conditions $\phi_0$ in $\V_N^c$.

We consider a sequence of i.i.d. real random variables
$(X_i)_{i\in{\mathbb N}}$, such that $X_i$ has probability density
$\exp(-V)dr$ on ${\mathbb R}$. Then $q = {\mathbb
E}\left[X_1^2\right]$, see \eqref{qdet}. For $n\in{\mathbb N}$ we
set $S_n:=X_1+\cdots+X_n$, $S_0:=0$. Moreover, for any $c\in\R$ and
$N\in\N$ we set
\[
T^{N,c}_i:= S_{i-1} - \frac1N \sum_{j=1}^{N-1} S_j + c\, N^{1/2},
\qquad i=1,\ldots,N,
\]
and
\[
\V_N^c \, := \, \left\{\phi\in\R^N:\sum_{i=1}^N \phi_i\, =\, c\,
N^{3/2}\right\} \, = \, \V_N+c\, N^{1/2}\, {\bf 1}.
\]
Notice that a.s. $T^{N,c}=(T^{N,c}_1,\ldots,T^{N,c}_N)\in\V_N^c$.
Clearly $\V_N^c$ is a $(N-1)$-dimensional affine subspace of $\R^N$;
we denote by $\cL^{N-1}(d\phi)$ the induced $(N-1)$-dimensional
Lebesgue measure.
\begin{lemma}\label{proba}
The law of $(T_1^{N,c},\ldots,T_N^{N,c})$ on $\V_N^c$ is
\begin{equation}\label{pn}
{\bf P}_N^c(d\phi) \, := \, \frac 1{Z_N^c} \, 1_{(\phi\in \V_N^c)}\,
\exp\left\{ - \cH_N(\phi) \right\} \, \cL^{N-1}(d\phi),
\end{equation}
where $Z_N^c$ is a normalization constant and $\cH_N$ is the
Hamiltonian
\[
\cH_N(\phi) \, := \, \sum_{x=2}^N V(\phi(x)-\phi(x-1)), \qquad
\phi\in\R^N.
\]
\end{lemma}
\noindent{\it Proof.} It is enough to prove the case $c=0$. We set
$\tau:\R^{n-1}\mapsto\R^N$,
\[
\tau(y) := -\frac 1N\sum_{k=1}^{N-1}y_k \cdot {\bf 1} +\left(0, y_1,
\ldots, y_{N-1} \right), \quad y\in\R^{N-1}.
\]
For all $f\in C_b(\R^N)$, we have
\begin{align*}
\E[f(T^{N,0})] = \int_{\R^{N-1}} f(\tau(y)) \,
e^{-V(y_1)-V(y_2-y_1)-\cdots-V(y_{N-1}-y_{N-2})} \, dy_1\cdots
dy_{N-1}.
\end{align*}
Now we define the $(N-1)\times(N-1)$ matrix
\[
L := (L_{ij}), \qquad L_{ij} = 1_{(i=j)} - \frac 1N,
\]
so that $\tau_i(y)=(Ly)_{i-1}$ for all $i=2,\ldots,N$. Let us now
use the following change of variable
\[
\R^{N-1} \ni y \mapsto (\phi_2,\ldots,\phi_N)\in \R^{N-1}, \qquad
\phi_i:= (Ly)_{i-1}, \quad i=2,\ldots,N.
\]
Moreover we set
\[
\phi_1:=-\frac 1N\sum_{k=1}^{N-1}y_k=-(\phi_2+\cdots+\phi_N).
\]
Then $(\phi_1,\ldots,\phi_N)\in \V_N$ and $y_1=\phi_2-\phi_1$,
$y_i-y_{i-1}=\phi_{i+1}-\phi_i$, for all $i=1,\ldots,N-1$. Finally
\[
\E[f(T^{N,0})] = \frac 1{|\det L|} \int_{\R^{N-1}}
f(\phi_1,\ldots,\phi_N) \,
e^{-V(\phi_2-\phi_1)-\cdots-V(\phi_N-\phi_{N-1})} \, d\phi_2\cdots
d\phi_N. \qed
\]

\medskip\noindent
We also set ${\bf P}_N^{c,+} = {\bf P}_N^c(\,\cdot\, | \,
\Omega_N^+)$. Then
\begin{equation}\label{pnc}
{\bf P}_N^{c,+}(d\phi) \, = \, \frac1{Z^{c,+}_N} \, 1_{(\phi\in
\V_N^c\cap\Omega_N^+)}\, \exp\left\{ - \cH_N(\phi) \right\} \,
\cL^{N-1}(d\phi),
\end{equation}
where $Z_N^{c,+}={\bf P}_N^c(\Omega_N^+)$ is a normalization
constant.

\medskip\noindent
Since $\V_N^c=c{\bf 1}+\V_N$ is an affine space obtained by a
translation of $\V_N$, it is natural to consider $\V_N$ as its
tangent space. More precisely, for any $F: \V_N^c\mapsto\R$ in
$C^1$, one can define a gradient $\nabla_{\V_N}F:
\V_N^c\mapsto \V_N$ as follows
\[
\left.\frac d{d\varepsilon} \, F(\phi+\gep \, v) \right|_{\gep=0} \,
= \, \langle\nabla_{\V_N}F(\phi),v\rangle_{\V_N}, \qquad \forall \
\phi\in\V_N^c, \ v\in\V_N,
\]
recall \eqref{gradiente}. Notice that $\nabla_{\V_N}$ is the
gradient operator in $\V_N$ with respect to the scalar product
$\langle\cdot,\cdot\rangle_{\V_N}$. If $F\in C^1(\R^N)$ and
$\phi\in\V_N^c$, then it is possible to compare the gradient in
$\V_N$ and the standard gradient $\nabla F=(\frac{\partial
F}{\partial \phi_i}, i=1,\ldots,N)$
\[
\nabla_{\V_N}F = \sigma\sigma^T\nabla F, \qquad \|\nabla_{\V_N}
F\|^2_{\V_N} = \|\sigma^T\nabla F\|^2_{\R^N} = \langle\nabla F,
\sigma\sigma^T\nabla F\rangle_{\R^N}.
\]

\begin{prop}\label{dirf} Let $c>0$.
\begin{enumerate}
\item The Markov process $(\phi(t,\phi_0))_{t\geq
0,\phi_0\in\V_N^c\cap\Omega_N^+}$ is the diffusion generated by the
symmetric Dirichlet Form in $L^2(\Omega_N^+,{\bf P}_N^{c,+})$,
closure of
\begin{align*}
C^1_b(\Omega_N^+) \, \ni \, F \, \mapsto \, e^{c,N}(F,F) \, & := \,
\int \sum_{x,y\in \Gamma_N} \frac{\partial F}{\partial \phi(x)} \,
[\sigma\sigma^T]_{xy} \, \frac{\partial F}{\partial \phi(y)} \,
d{\bf P}_N^{c,+} \\ & = \int \|\nabla_{\V_N}F\|^2_{\V_N} \, d{\bf
P}_N^{c,+}.
\end{align*}
\item ${\bf P}_N^{c,+}$ is the only tempered invariant probability measure of
$\phi$ on $\V_N^c\cap\Omega_N^+$, where temperedness means having finite second moment.
\end{enumerate}
\end{prop}
\noindent{\bf Proof}. Closability of $e^{c,N}$ on
$C^1_b(\Omega_N^+)$ follows from Theorem \ref{main1}, since the
Hamiltonian $\cH_N$ and the set $\V_N^c\cap\Omega_N^+$ are convex
and ${\bf P}_N^{c,+}$ is therefore log-concave (see Theorem 9.4.11
of \cite{ags}). Since $\V_N^c\cap\Omega_N^+$ is locally compact, by
Fukushima's theory of Dirichlet forms there exists a continuous
Markov process $(\psi_t,t\geq 0)$ in $\V_N^c\cap\Omega_N^+$,
starting from quasi-every $\psi_0\in\V_N^c\cap\Omega_N^+$, weak
solution of \eqref{sde}. By the pathwise uniqueness result of Lemma
\ref{exun}, $(\psi_t,t\geq 0)$ and $(\phi_t,t\geq 0)$ are identical
in law if $\psi_0=\phi_0$ and therefore $(\phi_t,t\geq 0,
\phi_0\in\V_N^c\cap\Omega_N^+)$ is the Markov process associated
with $e^{c,N}$.

The second assertion follows from point (e) of Theorem \ref{main1},
since ${\bf P}_N^{c,+}\in\ProbabilitiesTwo{\R^N}$ by the convexity of
$V$ and in particular \eqref{qdet}.
\qed

\section{The rescaling}\label{five}

Recall now the rescaling map $\Lambda_N:\R^N\mapsto L^2(0,1)$,
defined in \eqref{rescaling}. In this section we show how the scalar
product of $\V_N$ is transformed under this map. This issue is crucial
for the proof of \eqref{bka0} and \eqref{sig0} in our setting, see Proposition
\ref{asz2} below.

We define the linear subspace $H_N$ of $L^2(0,1)$ as the image of
$\Lambda_N$. We denote by $1_{I(x)}$ the indicator function of the
interval $I(x)$, where
\[
I(0) \, := \, \emptyset, \qquad I(x) \, := \, [(x-1)/N,x/N), \quad
x\in\Gamma_N.
\]
Then, by the definition of $\Lambda_N$
\[
H_N = \left\{ \sum_{i=1}^N a_i \, 1_{I_i}, \quad
(a_1,\ldots,a_N)\in\R^N\right\},
\]
i.e. $H_N$ can be identified with the space of functions on $[0,1)$
being constant on $I(x)$ for all $x\in\Gamma_N$.

Let $B$ denote a standard Brownian motion in $\R$ with $B_0=0$. We
set
\[
\overline B_N \, := \, \frac{B_{\frac1N}+B_{\frac2N}+\cdots+B_1}N, \qquad
\overline B \, := \, \int_0^1 B_r \, dr.
\]
Then we define the process
\[
Y^N_r \, := \, B_{\lfloor Nr+1\rfloor/N}
- \overline B_N, \qquad r\in[0,1),
\]
\[
Y_r \, := \, B_r-\overline B, \qquad r\in[0,1],
\]
where $\lfloor\cdot\rfloor$ denotes the integer part.
Notice that almost surely
\[
\langle Y^N,1\rangle \, = \, \langle Y,1\rangle \, = \, 0,
\qquad Y^N_r \, \to \, Y_r, \quad \forall \ r\in[0,1)
\]
as $N\to\infty$. Both processes are centered Gaussian. Recall that
$\langle\cdot,\cdot\rangle=\langle\cdot,\cdot\rangle_{L^2(0,1)}$
denotes the scalar product in $L^2(0,1)$. Now we define
\[
\langle h,k\rangle_{H_N} \, := \, \E\left[\langle h,Y^N\rangle \,
\langle k,Y^N\rangle \right] + \langle h,1\rangle \, \langle
k,1\rangle, \qquad \forall \ h,k\in H_N,
\]
\[
\langle h,k\rangle_H \, := \, \E\left[\langle h,Y\rangle\, \langle
k,Y\rangle \right] + \langle h,1\rangle \, \langle k,1\rangle,
\qquad \forall \ h,k\in L^2(0,1).
\]
\begin{lemma}\label{h1}$ $
\begin{itemize}
\item For any $N\in\N$ and $h\in H_N$
\begin{align}\label{equiv}
\langle h,h\rangle_{H_N} & \, = \, \langle h,1\rangle^2 + \frac 1N
\sum_{i=1}^{N-1} \left( \sum_{j=1}^i\langle h-\langle
h,1\rangle,1_{I(j)}\rangle\right)^2
\\ \nonumber & = \, \langle h,1\rangle^2 + \E\left[
\langle h,\Lambda_N D \rangle^2\right],
\end{align}
where $D$ is defined in \eqref{D}. In particular, if $h\ne 0$ then
$\langle h,h\rangle_{H_N}>0$.
\item For any $h\in L^2(0,1)$
\[
\langle h,h\rangle_H \, = \, \langle h,1\rangle^2 + \int_0^1 \left(
-\langle h,1\rangle+\int_0^t h(s) \, ds\right)^2 dt.
\]
In particular, if $h\ne 0$, then $\langle h,h\rangle_H>0$.
\end{itemize}
\end{lemma}
\noindent{\it Proof}. Let $h\in H_N$ and set
\[
k \, := \, \sum_i \langle h-\langle
h,1\rangle,1_{I(1)}+\cdots+1_{I(i)}\rangle \, 1_{I(i)},
\]
and notice that $\langle k,1_{I(N)}\rangle=0$. Then
\begin{align*}
& \langle h,h\rangle_{H_N} - \langle h,1\rangle^2 \, = \,
\E\left[\langle h-\langle h,1\rangle,B_{\lfloor N\cdot+1\rfloor/N}
\rangle^2 \right] \, = \, \E\left[\left(\sum_{i=1}^N \langle h-
\langle h,1\rangle,1_{(i)} \rangle B_{\frac iN}\right)^2 \right]
\\ &  = \,
\E\left[\left(\langle k,1_{I(N)}\rangle \, B_1 - \sum_{i=1}^{N-1}
\langle k,1_{I(i)}\rangle\left( B_{\frac{i+1}N}-B_{\frac
iN}\right)\right)^2 \right] \, = \, \frac 1N \sum_{i=1}^{N-1}
\langle k,1_{I(i)}\rangle^2,
\end{align*}
and \eqref{equiv} is proven, also recalling Lemma \ref{repre}.

Analogously, for any $h\in L^2(0,1)$ we set $k_r:=\int_0^r
(h-\langle h,1\rangle)$. Then we find $k_1=0$ and
\[
\langle h,h\rangle_H  - \langle h,1\rangle^2\, = \, \E\left[\langle
h-\langle h,1\rangle,B \rangle^2 \right] \, = \, \E\left[\left(k_1
\, B_1 - \int_0^1 k \ dB\right)^2 \right] \, = \, \int_0^1 k^2. \qed
\]
Therefore $\langle\cdot,\cdot\rangle_{H_N}$, respectively
$\langle\cdot,\cdot\rangle_H$, \ defines a scalar product on $H_N$,
resp. on $L^2(0,1)$. We define the Hilbert space $H$, completion of
$L^2(0,1)$ with respect to the scalar product
$\langle\cdot,\cdot\rangle_H$. Notice that the associated norms are
controlled by the $L^2(0,1)$ norm.
\begin{lemma}\label{control}
For all $N\in\N$ and $h\in H_N$
\[
\|h\|^2_{H_N} \leq \|h\|^2_{L^2(0,1)}.
\]
For all $h\in L^2(0,1)$
\[
\|h\|^2_{H} \leq  \|h\|^2_{L^2(0,1)}.
\]
\end{lemma}
\noindent{\it Proof}. For any $N\in\N$ and $h\in H_N$
\begin{align*}
& \langle h,h\rangle_{H_N} - \langle h,1\rangle^2 = \E\left[\langle
h-\langle h,1\rangle,B_{\lfloor N\cdot+1\rfloor/N} \rangle^2 \right]
\\ & \leq \|h-\langle h,1\rangle\|^2_{L^2(0,1)} \ \E\left[\|B_{\lfloor
N\cdot+1\rfloor/N} \|^2_{L^2(0,1)} \right] = \|h-\langle
h,1\rangle\|^2_{L^2(0,1)} \ \frac 1N \sum_{i=1}^N \frac iN \\  &
\leq \|h-\langle h,1\rangle\|^2_{L^2(0,1)}.
\end{align*}
Therefore
\[
\langle h,h\rangle_{H_N} \leq \langle h,1\rangle^2+ \|h-\langle
h,1\rangle\|^2_{L^2(0,1)} = \|h\|^2_{L^2(0,1)}.
\]
Analogously, for any $h\in L^2(0,1)$
\begin{align*}
& \langle h,h\rangle_{H} - \langle h,1\rangle^2 = \E\left[\langle
h-\langle h,1\rangle,B \rangle^2 \right] \leq \|h-\langle
h,1\rangle\|^2_{L^2(0,1)} \ \E\left[\|B \|^2_{L^2(0,1)} \right] \\ &
= \|h-\langle h,1\rangle\|^2_{L^2(0,1)} \
\int_0^1 t\, dt 
\leq \|h-\langle h,1\rangle\|^2_{L^2(0,1)}.
\qed
\end{align*}

\medskip\noindent We define now the image measures of
${\bf P}_N^c$ and ${\bf P}_N^{c,+}$ under $\Lambda_N$,
\[
\nu_N^c \, := \, \Lambda_N^* ({\bf P}_N^c), \qquad \nu_N^{c,+} \, :=
\, \Lambda_N^* ({\bf P}_N^{c,+}), \qquad c>0,
\]
where $\Lambda_N$, ${\bf P}_N^c$ and ${\bf P}_N^{c,+}$ are defined,
respectively, in \eqref{rescaling}, \eqref{pn} and \eqref{pnc}.
Finally, we set for all $c\in\R$
\[
H_N^c:=\left\{ h\in H_N, \ \langle h,1\rangle =c \right\}, \qquad
H^c:=\left\{ h\in H, \ \langle h,1\rangle =c \right\};
\]
in particular, $H_N^0$ and $H^0$ are Hilbert space w.r.t. to the
restrictions of $\langle\cdot,\cdot\rangle_{H_N}$, respectively
$\langle\cdot,\cdot\rangle_H$, that we denote
\[
\langle h,k\rangle_{H_N^0} \, := \, \E\left[\langle h,Y^N\rangle \,
\langle k,Y^N\rangle \right], \qquad \forall \ h,k\in H_N^0,
\]
\[
\langle h,k\rangle_{H^0} \, := \, \E\left[\langle h,Y\rangle\,
\langle k,Y\rangle \right], \qquad \forall \ h,k\in H^0.
\]
By \eqref{equiv} and Lemma \ref{repre}, we see that the scalar
product in $H^0_N$ is the push-forward of the scalar product in
$\V_N$ under $\Lambda_N$, i.e. for all $h\in H^0_N$
\begin{equation}\label{pushf}
\|h\|^2_{H_N^0} \, = \, \|\Lambda_N^{-1}h\|^2_{\V_N}.
\end{equation}

\medskip\noindent
As in the case of $\V_N^c$, for a differentiable $F:H_N^c\mapsto\R$
we can define a gradient $\nabla_{H_N^0}F:H_N^c\mapsto H_N^0$
\[
\left.\frac d{d\varepsilon} \, F(k+\gep \, h) \right|_{\gep=0} \, =
\, \langle\nabla_{H_N^0}F(k),h\rangle_{H_N^0}, \qquad \forall \ k\in
H_N^c, \ h\in H_N^0.
\]
Analogously for a differentiable $F:H^c\mapsto\R$ we can define a
gradient $\nabla_{H^0}F:H^c\mapsto H^0$
\[
\left.\frac d{d\varepsilon} \, F(k+\gep \, h) \right|_{\gep=0} \, =
\, \langle\nabla_{H^0}F(k),h\rangle_{H^0}, \qquad \forall \ k\in
H^c, \ h\in H^0.
\]
Moreover, $\Lambda_N:\V_N^c\mapsto H_N^c$ is bijective. Then, for
any $f\in C^1_b(H_N^c)$ we have $f\circ\Lambda_N\in C^1_b(\V^c_N)$
and
\begin{equation}\label{diff}
\sum_{x,y\in \Gamma_N} \frac{\partial (f\circ\Lambda_N)}{\partial
\phi(x)} \ [\sigma\sigma^T]_{xy} \ \frac{\partial
(f\circ\Lambda_N)}{\partial \phi(y)} \, = \, \frac 1{N^4} \,
\|\nabla_{H_N^0} f\|^2_{H_N^0}\circ \Lambda_N.
\end{equation}
Then we have for any $\varphi,\psi\in C^1_b(H_N^c)$
\[
{\cE}^{c,N}(f,g) \, := \, \int_{K_N} \langle\nabla_{H_N^0}
\varphi,\nabla_{H_N^0}\psi\rangle_{H_N^0} \, d\nu_N^{c,+} \, = \,
N^4 \, e^{c,N}(\varphi\circ\Lambda_N, \psi\circ\Lambda_N).
\]
We obtain readily from Proposition \ref{dirf}
\begin{prop}\label{h}
The bilinear form $({\cE}^{c,N},C^1_b(H_N^c))$ is closable in
$L^2(\nu_N^{c,+})$ and the closure $({\cE}^{c,N},D({\cE}^{c,N}))$ is
a symmetric Dirichlet form with associated Markov process
$\Phi^N$.
\end{prop}

\section{Proof of Theorem \ref{main}}

\subsection{The limit equation}

We recall that $B$ denotes a standard real Brownian motion and
\[
\overline B \, := \, \int_0^1 B_r \, dr,
\]
We define the process
\[
Y^c_\theta := q^{1/2}\left(B_\theta-\overline B\right)+c, \qquad
\theta\in[0,1],
\]
and $\nu^{c,+}$ as the law of $Y^c$ conditioned
to be non-negative on $[0,1]$. In other words, if $\nu^c$ is the law
of $Y^c$ and $K:=\{h\in L^2(0,1), h\geq 0\}$, then
$\nu^{c,+}=\nu^c(\,\cdot\, | K)$. The following result has been
proven in \cite{deza}.
\begin{prop}\label{deza} $ $ {\rm
\begin{enumerate}
\item For all $u_0\in H^c\cap K$
there exists a unique strong solution of \eqref{0}. We denote
$X_t(u_0):=u(t,\cdot)\in H^c\cap K$
\item The process $(X_t(u_0))_{t\geq 0,u_0\in H^c\cap K}$ is the diffusion associated with
the Dirichlet form $(\cE^c,D(\cE^c))$, closure of the symmetric
form
\[
{\cE}^c(\varphi,\psi) \, :=  \, \int \langle\nabla_{H^0} \varphi ,
\nabla_{H^0} \psi\rangle_{H^0} \, d\nu^{c,+}, \qquad \forall \
\varphi,\psi\in C^1_b(H^c).
\]
\item $\nu^{c,+}$ is the only invariant measure of $(X_t(u_0))_{t\geq 0,u_0\in H^c\cap K}$.
\end{enumerate}
}
\end{prop}

\subsection{Proof of \eqref{bka0} and \eqref{sig0}}
We are going to show now that, as $N\to\infty$, $\nu^{c,+}_N$
converges weakly to $\nu^{c,+}$ and the norm $\|\cdot\|_{H^0_N}$
converges to $\|\cdot\|_{H^0}$, in the sense of \eqref{bka0} and \eqref{sig0}.
\begin{prop}\label{asz2} In the notation of section \ref{five}
\begin{enumerate}
\item If $c>0$ then $\nu^{c,+}_N$ converges weakly in $H$ to $\nu^{c,+}$ as
$N\to+\infty$.
\item We have
\begin{equation}\label{bka}
\frac{1}{6} \, \Vert h\Vert_{H^0}\leq \Vert h\Vert_{H_N^0}\leq \Vert
h\Vert_{H^0} \qquad\forall \ h\in H_N^0, \ N \in\N.
\end{equation}
\item Denoting by $\Pi_N:H^0\to H_N^0$ the orthogonal
projections induced by the scalar product of $H^0$, we have
\begin{equation}\label{sig}
\lim_{N\to\infty}\Vert\Pi_N h\Vert_{H_N^0}=\Vert h\Vert_{H^0}
\qquad\forall \ h\in H^0.
\end{equation}
\end{enumerate}
\end{prop}
\noindent{\it Proof}. We start with weak convergence of
$\nu^{c,+}_N$ to $\nu^{c,+}$. We set $\nu_N^c \, := \, \Lambda_N^*
({\bf P}_N^c)$, i.e. $\nu_N^c$ is the law of the process $Y^{c,N}$
\[
Y^{c,N}_\theta \, := \, \frac{S_{\lfloor N\theta\rfloor}-\overline
S_N}{\sqrt N} \, + \, c, \qquad \theta\in[0,1).
\]
By the invariance principle, $\nu_N^c$ converges weakly to the law
$\nu^c$ of $Y^c:=q^{1/2}\left(B-\overline B\right)+c$, where $q$ is
defined in \eqref{qdet}. We have to prove now that for $c>0$
\[
\nu^c(\partial K) = \bbP\left(\inf_{\theta\in[0,1]} Y_\theta^c=0
\right) = 0.
\]
Notice that, by the symmetry of $Y^c$ with respect to time inversion
$\theta\mapsto1-\theta$, we have
\[
\bbP\left(\inf_{\theta\in[0,1]} Y_\theta^c=0 \right) \leq 2\,
\bbP\left(\inf_{\theta\in[0,1/2]} Y_\theta^c=0 \right).
\]
Notice that $\overline B\sim\cN(0,1/3)$.
By a standard Gaussian computation, it is easy to see that the law
of $(Y_\theta^c, \theta\in[0,1/2])$ is equivalent to the law of
\[
V_\theta:=q^{1/2}(B_\theta-Z) + c, \qquad \theta\in[0,1/2],
\]
where $Z\sim \cN(0,1/3)$ is independent of $B$. Since the minimum
value of $B$ over $[0,1/2]$ has the law of $|B_{1/2}|$, we obtain
that
\[
\bbP\left(\inf_{\theta\in[0,1/2]} V_\theta=0 \right)=
\bbP\left(|B_{1/2}|=Z-q^{-1/2}c\right)=0
\]
and therefore $\bbP\left(\inf_{\theta\in[0,1/2]} Y_\theta^c=0
\right)=0$. Then $\nu^c(\partial K) = 0$ and $\nu_N^c(\, \cdot \, |
K)=\nu^{c,+}_N$ converges weakly to $\nu^c(\, \cdot \, |
K)=\nu^{c,+}$.

\medskip\noindent
We prove now \eqref{bka} and \eqref{sig}. The key result is the
following lemma.
\begin{lemma}\label{h0}
For all $N\in\N$ and $h\in H_N$
\begin{equation}\label{est0}
\|h\|_{H_N}^2 + \frac1{6N^2}
\langle h,1\rangle^2 \, = \,
\|h\|_{H}^2 + \frac1{6N^2} \|h\|^2_{L^2(0,1)}.
\end{equation}
\end{lemma}
\noindent{\it Proof}. Since $\langle h,1\rangle_H= \langle
h,1\rangle_{H_N}=\langle h,1\rangle$, then (\ref{est0}) is
equivalent to
\[
\|h-\langle h,1\rangle\|_{H_N}^2 \, = \,
\|h-\langle h,1\rangle\|_{H}^2 + \frac1{6N^2} \|h-\langle h,1\rangle\|^2_{L^2(0,1)},
\qquad \forall \ h\in H_N.
\]
This, in turn, is equivalent to
\[
\E\left[\langle h,B_{\lfloor N\cdot+1\rfloor/N} \rangle^2 \right] \,
= \, \E\left[\langle h,B \rangle^2 \right] + \frac1{6N^2}
\|h\|^2_{L^2(0,1)}, \qquad \forall \ h\in H_N^0.
\]
This formula can be proven by noting that for all $i=1,\ldots,N$
\[
B_{\frac iN} \, = \, N\int_{\frac{i-1}N}^{\frac iN} B_s \, ds +
N\int_{\frac{i-1}N}^{\frac iN} \left(B_{\frac iN}-B_s\right) \, ds.
\]
Indeed, it follows that for all $h\in H_N$
\begin{align*}
& \E\left[\langle h,B_{\lfloor N\cdot+1\rfloor/N} \rangle^2 \right]
\, = \, \E\left[\left(\sum_{i=1}^N \langle h,1_{(i)} \rangle
B_{\frac iN}\right)^2 \right]
\\ &  = \,
\E\left[\langle h,B \rangle^2 \right] + \,
\E\left[\left(\sum_{i=1}^N \langle h,1_{(i)} \rangle
N\int_{\frac{i-1}N}^{\frac iN} \left(B_{\frac iN}-B_s\right) \,
ds\right)^2 \right]
\\ & \quad + \, 2 \, N^2 \, \E\left[\langle h,B \rangle\sum_{i,j=1}^N \langle h,1_{(i)} \rangle \,
\langle h,1_{(j)} \rangle \int_{\frac{j-1}N}^{\frac jN} B_r \, dr
\int_{\frac{i-1}N}^{\frac iN} \left(B_{\frac iN}-B_s\right) ds\right]
\end{align*}
By independence of increments of the Brownian motion, the
second term in the right hand side is
\[
\E\left[\left(\sum_{i=1}^N \langle h,1_{(i)} \rangle
N\int_{\frac{i-1}N}^{\frac iN} \left(B_{\frac iN}-B_s\right) \,
ds\right)^2 \right] = \, \frac1{3N} \sum_{i=1}^N \langle h,1_{(i)}
\rangle^2 \, = \, \frac1{3N^2} \|h\|^2_{L^2(0,1)}.
\]
Now, for the third term, we need to
calculate
\[
I_{ij} \, := \, \E\left[\int_{\frac{j-1}N}^{\frac jN} B_r \, dr
\int_{\frac{i-1}N}^{\frac iN} \left(B_{\frac iN}-B_s\right)
ds\right].
\]
Again by independence we have $I_{ij}=0$ if $j<i$. On the other hand
\[
i<j \, \Longrightarrow \,
I_{ij} \, = \, \int_{\frac{j-1}N}^{\frac jN} dr
\int_{\frac{i-1}N}^{\frac iN} \left(\frac iN-s\right) ds \, = \, \frac1{2N^3},
\]
\[
i=j \, \Longrightarrow \,
I_{ii} \, = \, \int_{\frac{i-1}N}^{\frac iN} dr
\int_{\frac{i-1}N}^{\frac iN} \left(s-r\right) ds \, = \, \frac1{6N^3}.
\]
Then we must compute for all $h\in H_N$
\begin{align*}
& \frac 1N
\sum_{i<j} \langle h,1_{(i)} \rangle \, \langle h,1_{(j)} \rangle +
\frac1{3N} \sum_i \langle h,1_{(i)} \rangle^2 \, = \,
\frac 1{2N}\sum_{i\neq j} \langle h,1_{(i)} \rangle \, \langle h,1_{(j)} \rangle +
\frac1{3N} \sum_i \langle h,1_{(i)} \rangle^2
\\ & = \frac 1{2N}\sum_{i,j} \langle h,1_{(i)} \rangle \, \langle h,1_{(j)} \rangle
- \frac1{6N} \sum_i \langle h,1_{(i)} \rangle^2 \, = \,
\frac 1{2N}\langle h,1\rangle^2 - \frac1{6N^2} \|h\|^2_{L^2(0,1)}.
\end{align*}
Finally, we have proven that for all $h\in H_N$
\[
\E\left[\langle h,B_{\lfloor N\cdot+1\rfloor/N} \rangle^2 \right] \,
= \, \E\left[\langle h,B \rangle^2 \right] + \, \frac1{6N^2}
\|h\|^2_{L^2(0,1)} \, + \, \frac 1{2N}\langle h,1\rangle^2
\]
and choosing $h$ such that $\langle h,1\rangle=0$ we have the
desired result. \qed

\medskip\noindent
{\it End of the proof of Proposition \ref{asz2}}. We prove now
\eqref{bka}, namely the estimate
\begin{equation}\label{h00}
\frac16 \, \|h\|_{H_N^0}^2 \,  \leq \, \|h\|_{H^0}^2 \, \leq \,
\|h\|_{H_N^0}^2, \qquad \forall \ N\in\N, \ h\in H_N^0.
\end{equation}
The second inequality of \eqref{h00} follows from (\ref{est0}). For
the first inequality, recall now (\ref{equiv}), where we proved that
for all $h\in H_N^0$
\[
\|h\|^2_{H_N^0} \, = \, \frac 1N\sum_{i=1}^{N-1} \left(\sum_{j=1}^i
\langle 1_{(j)},h\rangle\right)^2.
\]
Then we obtain for all $h\in H_N^0$
\begin{align*}
& \|h\|^2_{L^2(0,1)} \, = \, N\sum_{i=1}^N \langle
1_{(i)},h\rangle^2 = \, N\sum_{i=1}^{N-1} \left(\sum_{j=1}^i \langle
1_{(j)},h\rangle- \sum_{j=1}^{i-1} \langle 1_{(j)},h\rangle\right)^2
+ N \left(\sum_{j=1}^{N-1} \langle 1_{(j)},h\rangle\right)^2
\\ & \leq \, 4\, N\sum_{i=1}^{N-1}
\left(\sum_{j=1}^i \langle 1_{(j)},h\rangle\right)^2+ N
\left(\sum_{j=1}^{N-1} \langle 1_{(j)},h\rangle\right)^2\, \leq \, 5
\, N^2 \, \|h\|^2_{H_N^0}.
\end{align*}
Using (\ref{est0}) we obtain the first inequality and \eqref{h00} is
proven.

\medskip\noindent
We prove now \eqref{sig}, namely we prove that, denoting by
$\Pi_N:H^0\to H_N$ the orthogonal projections induced by the scalar
product of $H^0$, we have
\[
\lim_{N\to\infty}\Vert\Pi_N h\Vert_{H_N^0}=\Vert h\Vert_{H^0}
\qquad\forall \ h\in H^0.
\]
We denote by $P_N:L^2(0,1)\mapsto L^2(0,1)$ the following
projection
\begin{equation}\label{pi_N}
P_N h \, := \, \sum_{i=1}^N N \, \langle h,1_{I(i)}\rangle \,
1_{I(i)}, \qquad h\in L^2(0,1).
\end{equation}
Then $P_N$ is an orthogonal projector with respect to the scalar
product of $L^2(0,1)$ and for all $h\in L^2(0,1)$,
$\|h-P_Nh\|_{L^2(0,1)}\to 0$ as $N\to\infty$. Now, let us fix $h\in
L^2(0,1)\cap H^0$; then we have
\begin{equation}\label{u1}
\Vert P_N h\Vert_{H_N^0}^2 = \E\left[ \langle Y^N,h\rangle^2\right]
\to \E\left[ \langle Y,h\rangle^2\right] = \Vert h\Vert_{H^0}^2,
\qquad N\to\infty.
\end{equation}
Now we claim that $\|\Pi_N h\Vert_{H^0}^2 \to \Vert h\Vert_{H^0}^2$,
as $N\to\infty$. Indeed, $\Pi_N$ is the element of minimal
$H^0$-distance from $h$ in $H_N^0$. Then, since $P_Nh$ belongs to
$H_N^0$, by Lemma \ref{control}
\begin{equation}\label{u2}
\Vert \Pi_Nh-h\Vert_{H^0} \leq \Vert P_Nh-h\Vert_{H^0} \leq \Vert
P_Nh-h\Vert_{L^2(0,1)} \to 0, \qquad N\to\infty.
\end{equation}
Now, by \eqref{est0}
\[
\Vert \Pi_N h\Vert_{H_N^0}^2 = \Vert \Pi_N h\Vert_{H^0}^2 + \frac
1{6N^2}\|\Pi_N h\|^2_{L^2(0,1)} \geq \Vert \Pi_N h\Vert_{H^0}^2 \to
\Vert h\Vert_{H^0}^2, \qquad N\to\infty.
\]
In particular
\[
\liminf_{N\to\infty} \Vert \Pi_N h\Vert_{H_N^0} \geq \Vert
h\Vert_{H^0}.
\]
On the other hand, by \eqref{h00}
\[
\Vert \Pi_N h\Vert_{H_N^0} \leq \Vert P_N h\Vert_{H_N^0} + \Vert
P_N h -\Pi_N h\Vert_{H_N^0} \leq \Vert P_N h\Vert_{H_N^0} +
\Vert P_N h -\Pi_N h\Vert_{H^0}.
\]
Since $\lim_N (P_N h -\Pi_N h)=0$ in $H^0$ by \eqref{u2}, then by
\eqref{u1} we find
\[
\limsup_{N\to\infty} \Vert \Pi_N h\Vert_{H_N^0} \leq \Vert
h\Vert_{H^0}.
\]
If we set now
\[
\psi_N:H^0\mapsto\R, \qquad \psi_N(h)=\Vert \Pi_N h\Vert_{H_N^0},
\]
then $\psi_N$ is Lipschitz-continuous in the $H^0$-norm uniformly in $N$, since
\[
\Vert \Pi_N h\Vert_{H_N^0} \leq \Vert \Pi_N h\Vert_{H^0}  \leq \Vert h\Vert_{H^0}
\]
by \eqref{bka} and by the definition of $\Pi_N$.
Moreover and $\psi_N(h)\to\Vert h\Vert_{H^0}$ as $N\to\infty$
for all $h$ in $L^2(0,1)\cap H^0$. Since $L^2(0,1)\cap H^0$ is dense in $H^0$,
this concludes the proof of Proposition \ref{asz2}. \qed

\subsection{Proof of Theorem \ref{main}}

In order to prove Theorem \ref{main}, it is now enough
to notice that by Propositions \ref{h}, \ref{deza} and \ref{asz2},
Theorems \ref{main1} and \ref{main3} apply and yield the
desidered convergence result.

\end{document}